\documentclass[11pt]{amsart}
\usepackage{graphicx}

\begin{document}

\centerline {\Large \bf Fibonacci-Lucas densities}

\

\

\

\centerline {Mihai Caragiu} \centerline {Department of
Mathematics, Ohio Northern University} 
\centerline {Ada, OH 45810}
\centerline {m-caragiu1@onu.edu}

\

\centerline {Jacob L. Johanssen} 
\centerline {Department of
Mathematics, Ohio Northern University}
\centerline {Ada, OH 45810}
\centerline {j-johanssen@onu.edu}

\

\noindent {\bf ABSTRACT}. Both Fibonacci and Lucas numbers can be described combinatorially in terms of $0-1$ strings without consecutive ones. In the present article we explore the occupation numbers as well as the correlations between various positions in the corresponding configurations.

\

\noindent (2000) Mathematics Subject Classification: 11B39, 05A15

\noindent Keywords: Fibonacci numbers, Lucas numbers, occupation numbers.

\

\section {Introduction}

\

\noindent The Fibonacci numbers have found interesting applications in physical sciences [1],[2],[3],[4]. In the present paper we will start from the standard combinatorial interpretation of the Fibonacci sequence
$$F_0=0,\ \ F_1=1,\ \ F\sb n=F\sb {n-1}+F\sb {n-2}\ \ \ {\rm for \ n\geq 2}$$
in terms of binary sequences: $F\sb {n+2}$ represents the number of binary ($0-1$) sequences of length $n$ with no two 1's adjacent. Inspired by this, we define a {\it Fibonacci device} to be a linear arrangement of $n$ cells 
$C\sb 1,...,C\sb n$, with each cell $C\sb i$ hosting a local variable $\sigma \sb i\in \{ 0,1\}$. The set of admissible states of the Fibonacci device are the $n$-tuples
$$\sigma = \left (  \sigma \sb 1,\sigma \sb 2,...,\sigma \sb n \right )\in \{ 0,1\}\sp n\leqno (1)$$
subjected to the constraint
$$\sigma \sb i\sigma \sb {i+1}=0,\ \ \ {\rm for \ i=1,...,n-1}.\leqno (2)$$
Let us agree to call the $n$-tuples (1) subjected to the constraint (2), {\it Fibonacci states}. We will denote by $\Phi\sb n$ the set of Fibonacci states of a $n$-cell Fibonacci device. Thus $\left | \Phi \sb n \right |=F\sb {n+2}$.

\

\noindent The Fibonacci devices can be seen as physical (toy-)models for the bit-string interpretation of the Fibonacci numbers. The above specifications make possible the calculation of various statistical averages of functions of state
$$f:\Phi \sb n \to \mathbb R.\leqno (3)$$
The average $\langle f \rangle$ of (3) is defined by
$$\langle f \rangle=\frac {1}{F\sb {n+2}}\sum \sb {\sigma \in \Phi \sb n}f(\sigma).\leqno (4)$$
We will apply the general principle (4) to calculate various meaningful averages involving Fibonacci devices, such as the average occupation number of a given cell and the correlations between cells. A similar analysis will be performed for Lucas devices, which are similar to Fibonacci devices except that we assume `periodic boundary conditions' (the cells are wrapped around a circle).

\

\section {Fibonacci densities}

\noindent From (4) it follows that the average occupation number $\langle \sigma \sb i \rangle $ of a given cell $C\sb i$ of a Fibonacci device is given by 
$$\langle \sigma \sb i \rangle = \frac {1}{F\sb {n+2}}\sum \sb {\sigma \in \Phi \sb n}\sigma \sb i.\leqno (5)$$
Note that in (5) only the Fibonacci states corresponding to $\sigma \sb i=1$ (case in which $\sigma \sb {i-1}=\sigma \sb {i+1}=0$ necessarily) contribute to the sum. Therefore
$$\langle \sigma \sb i \rangle =\frac {1}{F\sb {n+2}} \sum \sb {\sigma \in \Phi \sb n, \ \sigma \sb i=1}\sigma \sb i. \leqno (6)$$
By using the combinatorial interpretation of the Fibonacci numbers specified in the introduction, we can evaluate the sum in (6) and rewrite the average occupation number $\langle \sigma \sb i \rangle$, which can also be seen as a local `Fibonacci density' $$\langle \sigma \sb i \rangle=\rho \sb{n,i}\sp F,$$ as
$$\rho \sb{n,i}\sp F =\frac {F\sb i  \cdot F_{n - i + 1}}{F\sb {n+2}}.\leqno (7)$$
By using Binet's formula
$$F\sb n=\frac {\alpha \sp n -\beta \sp n}{\sqrt 5},\leqno (8)$$
where $\alpha =\frac {1+\sqrt 5}{2}$ and $\beta =\frac {1-\sqrt 5}{2}$, we can rewrite the density function (8) as 
$$\rho \sb {n,i}\sp F  = \frac{1}
{{\sqrt 5 }}\frac{{\left( {\alpha ^i  - \beta \sp i } \right) \cdot \left( {\alpha \sp {n - i + 1}  - \beta \sp {n - i + 1} } \right)}}
{{\left( {\alpha \sp {n + 2}  - \beta \sp {n + 2} } \right)}}
.\leqno (9)$$
Note that if $i$ is fixed and $n$ grows larger and larger then the limit density at the $i$-th cell is
$$\rho \sb i\sp F  = \mathop {\lim }\limits_{n \to \infty } \rho \sb {n,i}\sp F  = \frac{{\alpha \sp i  - \beta \sp i }}
{{\sqrt 5 \alpha \sp {i + 1} }}
$$
Let us write the limit density $\rho \sb i\sp F$ in the form
$$\rho \sb i\sp F=\left( {\frac{1}
{{\alpha \sqrt 5 }}} \right)\left[ {1 - \left( {\frac{\beta }
{\alpha }} \right)\sp i } \right]\leqno (10)$$
Then (10), together with the fact that $$ - 1 < \frac{\beta }
{\alpha } < 0$$
shows an oscillating behavior of $\rho \sb i\sp F$. This explains the shape of the plot given in Figure 1, where we represented the local density function for a Fibonacci device with 80 cells (notice the oscillating behavior of the density when we approach the ends of the device).

\centerline {
\includegraphics[height=3in]{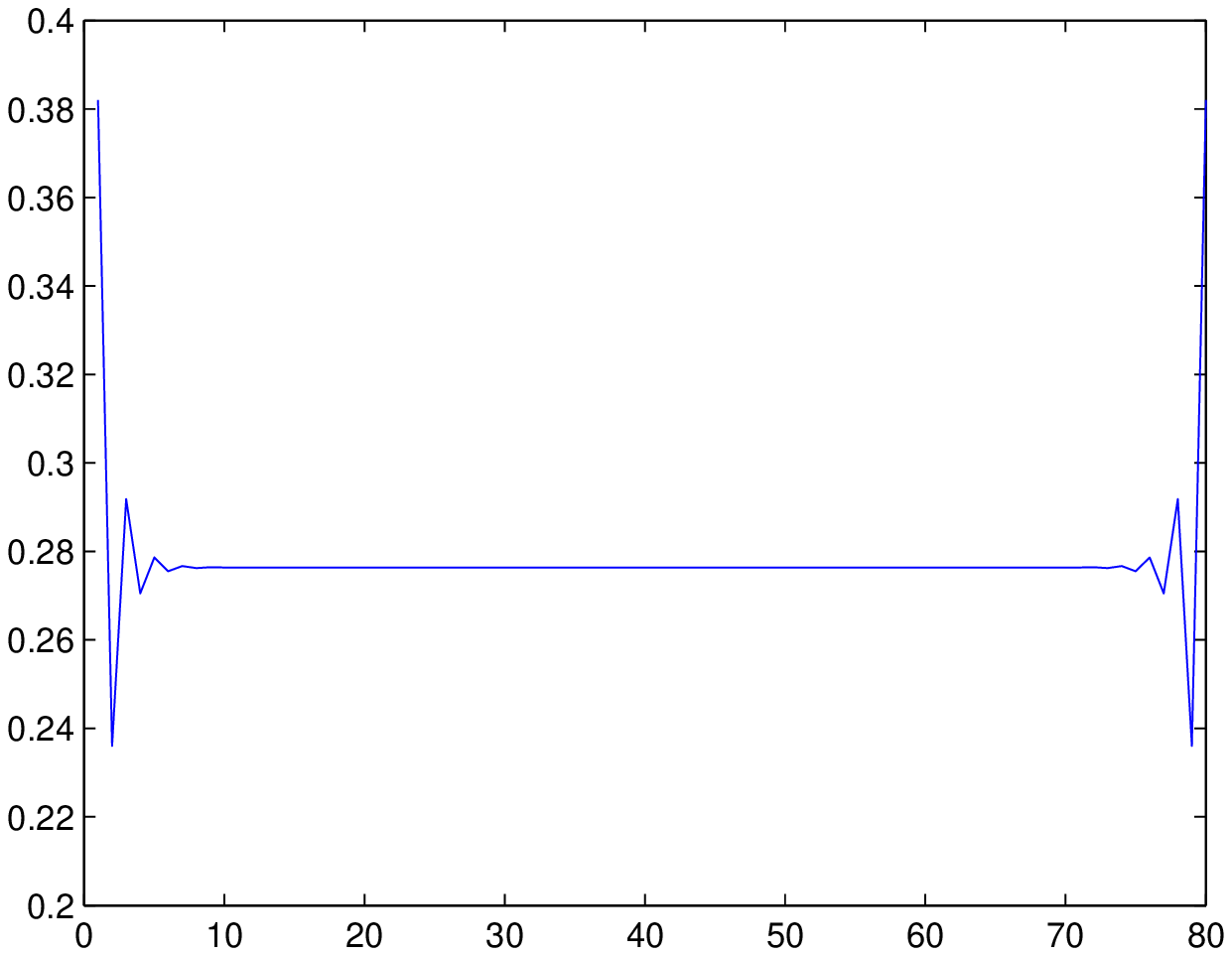}
}
\centerline {FIGURE 1. Local density $i\mapsto \rho \sb{80,i}\sp F$ for a Fibonacci device with 80 cells}

\

\noindent Note that deep inside the Fibonacci device, the limit density will be
$$\rho \sp F  = \frac{1}
{{\alpha \sqrt 5 }} = {\text{.2763932022}}...\leqno (11)
$$ 

\

\section {Lucas densities}

\noindent The Lucas numbers
$$L\sb 0=2,\ \ L\sb 1=1, \ \ L\sb n=L\sb {n-1}+L\sb {n-2}\ \ \ {\rm for \ n\geq 2}$$
have a combinatorial interpretation which is similar to that of the Fibonacci numbers: $L\sb {n}$ is the number of binary ($0-1$) sequences of length $n$ wrapped around a circle (so that the $n+1$-th position coincides with the with the first position, etc) with no two 1's adjacent. The analogue of the Binet's formula is
$$L\sb n=\alpha \sp n+\beta \sp n\leqno (12)$$
A a {\it Lucas device} is defined in a similar way with a Fibonacci device, the only difference being that the cells of a Lucas device are wrapped around a circle. Thus the {\it Lucas states} are the $n$-tuples
$$\sigma = \left (  \sigma \sb 1,\sigma \sb 2,...,\sigma \sb n \right )\in \{ 0,1\}\sp n$$
subjected to the constraint $\sigma \sb i\sigma \sb {i+1}=0$ for $i=1,...,n$, where $\sigma \sb {n+1}=\sigma \sb 1$.

\

\noindent If we denote by $\Lambda \sb n$ the set of Lucas states for a Lucas device with $n$ cells, the analogue of (4) will be
$$\langle f \rangle=\frac {1}{L\sb n}\sum \sb {\sigma \in \Lambda \sb n}f(\sigma),$$
where $\langle f \rangle$ is the statistical average of a state function $f:\Lambda \sb n \to \mathbb R$. The symmetry of a Lucas device implies that the average occupation number is the same for every cell $\sigma \sb i$ of a Lucas device with $n$ cells: $\rho \sb {n,i}\sp L=\rho \sb n \sp L$. If we fix a particular cell $C\sb i$ in a Lucas device with $n$ cells, there are exactly $F\sb {n-1}$ Lucas states for which $\sigma \sb i=1$ (which automatically implies that $\sigma \sb {i-1}=\sigma \sb {i+1}=0$). Therefore we have the following expression for the average occupation number (local density) for any cell in a Lucas device with $n$ cells:
$$
\rho \sb n \sp L=\frac {F\sb {n-1}}{L\sb n}.\leqno (13)
$$
From (12) and (13) we get
$$
\rho \sb n\sp L  = \frac{1}
{{\sqrt 5 }}\frac{{\alpha \sp {n - 1}  - \beta \sp {n - 1} }}
{{\alpha \sp n  + \beta \sp n }}
.
\leqno (14)$$
\noindent As a consequence of (14) we find that the limit density $\rho \sp L$ in the case of Lucas devices is the same as $\rho \sp F$:
$$\rho \sp L= \mathop {\lim }\limits_{n \to \infty } \rho \sb n\sp L  = \frac{1}
{{\alpha \sqrt 5 }} = .2763932022...$$

\

\section {Pair correlations in Fibonacci and Lucas devices}

\noindent Let $u,v \in \left\{ {0,1} \right\}$ be two bits. Define the correlation number
$\gamma \left( {u,v} \right)$ to be 1 if $u,v$ are identical and $-1$ if $u,v$ are distinct:
$$
\gamma \left( {u,v} \right) = \left\{ {\begin{array}{*{20}c}
   {1,{\text{ if }}u = v}  \\
   { - 1,{\text{ if }}u \ne v}  \\
 \end{array} } \right.
$$
This induces naturally, for every $1\leq k < l\leq n$, a state function 
$$\Gamma \sb {n,k,l}\sp F:\Phi \sb n \to \mathbb R$$
defined as follows:
$$\Gamma \sb {n,k,l}\sp F(\sigma )=\gamma (\sigma \sb k, \sigma \sb l).$$
Let us define the correlation function $C\sb {n,k,l}\sp F $ between the cell $C\sb k$ and the cell $C\sb l$ in a Fibonacci device with $n$ cells to be the statistical average of the state function $\Gamma \sb {n,k,l}\sp F$, that is, the average correlation number of the variables $\sigma \sb k$ and $\sigma \sb l$:
$$C\sb {n,k,l}\sp F =\langle \Gamma \sb {n,k,l}\sp F \rangle $$
From the definition of the correlation number and from the fact that there are $F\sb {n+2}$ Fibonacci states in a Fibonacci device with $n$ cells, we have:
$$C_{n,k,l}^F  = \frac{{N_{00}\sp F  + N_{11}\sp F  - N_{01}\sp F  - N_{10}\sp F }}
{{F_{n + 2} }},\leqno (15)$$
where for $i,j\in \{ 0,1\}$, $N\sb {ij}\sp F$ represents the number of Fibonacci states in a Fibonacci device with $n$ cells with $\sigma \sb k=i$ and $\sigma \sb l=j$. It is not hard to see that for $1\leq k<l\leq n$ we have $N\sb{00}\sp F  = F\sb{k + 1} F\sb {l - k + 1} F\sb {n - l + 2}$, $N\sb {11}\sp F  = F\sb k F\sb {l - k - 1} F\sb {n - l + 1}$, $N\sb {01}\sp F  = F\sb {k + 1} F\sb {l - k} F\sb {n - l + 1}$, and $N\sb {10}\sp F  = F\sb k F\sb {l - k} F\sb {n - l + 2}$. This, together with (15), implies
\scriptsize $$C\sb {n,k,l}\sp F  = \frac{{F\sb {k + 1} F\sb {l - k + 1} F\sb {n - l + 2} + F\sb k F\sb {l - k - 1} F\sb {n - l + 1}-F\sb {k + 1} F\sb {l - k} F\sb {n - l + 1}  - F\sb k F\sb {l - k} F\sb {n - l + 2} }}
{{F\sb {n + 2} }}. \leqno (16) $$
\normalsize
\noindent Note that if we use (16) to calculate $C\sb {n,k,l}\sp F$ for $k=l$ we will get the value to be expected, 1 (the self-correlation of any local variable). Indeed, by using $F\sb {-1}=1$ together with the Fibonacci identity $F\sb {a-1}F\sb b+F\sb aF\sb {b+1}=F\sb {a+b}$ ([5], p.9) for $a=k+1$ and $b=n-k+1$, we will obtain:
$$C\sb {n,k,k}\sp F=\frac{{F\sb {k + 1} F\sb {n - k + 2} + F\sb k F\sb {n - k + 1} }}
{{F\sb {n + 2} }}=\frac{F\sb {n+2}}{F\sb {n + 2}}=1.$$
Also, if we use (16) to calculate $C\sb {n,k,k+1}\sp F$ we will get
$$C\sb {n,k,k+1}\sp F=\frac{F\sb {k+1}F\sb {n-k+1}-F\sb {k+1}F\sb {n-k}-F\sb k F\sb {n-k+1}}{F\sb {n + 2} }=\frac{F\sb {k+1}F\sb {n-k-1}-F\sb k F\sb {n-k+1}}{F\sb {n + 2} }$$
By using $F\sb {k+1}=F\sb k+F\sb {k-1}$ in the last fraction above, we will find
$$C\sb {n,k,k+1}\sp F=-\frac{F\sb {k}F\sb {n-k}-F\sb {k-1} F\sb {n-k-1}}{F\sb {n + 2} }<0,$$
so that the correlation between two neighboring cells in the Fibonacci device is negative, as one should intuitively expect: for nearest-neighbor cells a negative correlation can be obtained in two ways (one of the two cells hosts a zero, with the other cell hosting a one) while a positive correlation can be obtained only in one way, namely if the local variables at both cells are zero.

\

\noindent One might also notice that $C\sb {n,k,l}\sp F$ given by (16) is invariant under the transformation $(k,l)\mapsto (n+1-l,n+1-k)$. This also follows from the following very intuitive reason: $C\sb {n,n+1-l,n+1-k}\sp F$ represents the correlation between the $k$-th and $l$-th cells {\it measured from the other end} of the Fibonacci device, so that $C\sb {n,k,l}\sp F=C\sb {n,n+1-l,n+1-k}\sp F$.

\

\noindent A similar correlation analysis can be performed for Lucas devices. In this case, if we exploit the circular symmetry of the Lucas device, it will be enough to calculate the correlation between cell 1 and cell $k$, where $1< k\leq \left\lceil {\frac{{n + 1}}
{2}} \right\rceil $. The appropriate state function will be 
$$\Gamma \sb {n,k}\sp L:\Lambda \sb n \to \mathbb R, \ \ \Gamma \sb {n,k}\sp L(\sigma )=\gamma (\sigma \sb 1, \sigma \sb k).$$
The correlation between cell 1 and cell $k$ in a Lucas device will be the average of $\Gamma \sb {n,k}\sp L$ and will be given by
$$C\sb {n,k}\sp L =\langle \Gamma \sb {n,k}\sp L \rangle = \frac{{N_{00}\sp L  + N_{11}\sp L  - N_{01}\sp L  - N_{10}\sp L }}
{{L\sb n }},\leqno (17)$$
where for $i,j\in \{ 0,1\}$, $N\sb {ij}\sp L$ represents the number of Lucas states in a Lucas device with $n$ cells with $\sigma \sb 1=i$ and $\sigma \sb k=j$. That is, $N\sb{00}\sp L  = F\sb{n-k + 2} F\sb {k}$, $N\sb {11}\sp L  = F\sb {n-k} F\sb {k-2}$, and $N\sb {01}\sp L  = N\sb {10}\sp L = F\sb {n-k + 1} F\sb {k-1}$. This, together with (17), gives us an expression for the correlation coefficient between the first and the $k$-th cell in a Lucas device with $n$ cells:
$$C\sb {n,k}\sp L  = \frac{{F\sb{n-k + 2} F\sb {k} + F\sb {n-k} F\sb {k-2}-2F\sb {n-k + 1} F\sb {k-1} }}
{{L\sb n }}. \leqno (18) $$
Mote that (18) works perfectly well in the case $k=1$ too. If $k=1$, $C\sb {n,1}\sp L$ represents the self-correlation of a cell, which we clearly expect to be 1. Indeed, if we set $k=1$ in (18), we get (after using $F\sb {-1}=1$ and $L\sb n=F\sb {n-1}+F\sb {n+1}$):
$$C\sb {n,1}\sp L  = \frac{F\sb {n+1}+F\sb {n-1}}{L\sb n }=1.$$
Also, $C\sb {n,2}\sp L$ represents the correlation between two nearest neighbors in a Lucas device with $n$ cells. By setting $k=2$ in (18) we get:
$$C\sb {n,2}\sp L  = \frac{F\sb n - 2F\sb {n-1}}{L\sb n }=\frac{F\sb {n-2} - F\sb {n-1}}{L\sb n }<0.$$
As in the case of the Fibonacci devices, this is consistent to the intuitive perception of a negative correlation between nearest neighbors. However if we set $k=3$ in (18) we find a {\it positive} correlation between two cells separated by another cell of a Lucas device with $n$ cells:
$$C\sb {n,3}\sp L  = \frac{2F\sb {n-1} +F\sb {n-3}-2F\sb {n-2}}{L\sb n }=\frac{3F\sb {n-3}}{L\sb n }>0.$$
 
\newpage

\

\noindent {REFERENCES}

\

\noindent [1] D'Amico, A., Faccio, M. and Ferri, G., {\it Ladder network characterization and Fibonacci numbers}, Nuovo Cimento D (1), 12 (1990), 1165--1173.

\

\noindent [2] Hoggatt, V. E., Jr. and Bicknell-Johnson, Marjorie, {\it Reflections across two and three glass plates}, Fibonacci Quarterly 17 (1979), no. 2, 118--142.

\

\noindent [3] Gumbs, Godfrey and Ali, M. K., {\it Electronic properties of the tight-binding Fibonacci Hamiltonian}, J. Phys. A 22 (1989), no. 8, 951--970.

\

\noindent [4] Tracy, Craig A., {\it Universality class of a Fibonacci Ising model}, J. Statistical Physics 51 (1988), no. 3-4, 481--490.

\

\noindent [5] N. N. Vorbobiev, {\it Fibonacci Numbers} (translated from the Russian by Mircea Martin), Birkhauser, 2002

\

\

\noindent Corresponding author address:

\

{\it

\noindent Mihai Caragiu

\noindent Department of Mathematics

\noindent Ohio Northern University

\noindent 262 Meyer Hall, Ada, OH 45810, USA

}

\end {document}